\documentclass{article}
\usepackage{amsmath}
\usepackage{graphicx}
\usepackage{enumitem}

\renewcommand{\vec}[1]{\boldsymbol{#1}}
\newcommand {\bJ}{{{J}_i}}
\newcommand {\bJT}{{{J}_i^T}}
\newcommand {\bJstar}{{{J}^{*}_i}}
\newcommand {\bJstarT}{{{J}^{*T}_i}}
\newcommand {\bK}{{{K}_i}}

\newcommand {\that}{\hat{\theta}}

\newcommand {\ei}{\epsilon_{i}}
\newcommand {\sm}{\sum_{i=1}^{n}}

\newcommand {\trgfi}{\left ( \frac{\bigtriangledown f}{f} \right
)_{i}}

\newtheorem{Res}{Result}

\begin{document}
\begin{center}
{\Large {\bf On comparison of  estimators for proportional error nonlinear regression models  in the limit of small measurement error}}

\vspace{0.5 cm}

{\large Richard A. Lockhart}$^*$

{\small{\it{Department of Statistics and Actuarial Science, Simon Fraser 
University, Burnaby, B.C. V5A 1S6, Canada}}}

\vspace{0.5 cm}
{\large Chandanie W. Navaratna}

{\small{\it{Department of Mathematics,
The Open University of Sri Lanka,
Nawala, Nugegoda,
Sri Lanka}}}

\end{center}

\begin{center}
Abstract
\end{center}

\noindent

In this paper, we compare maximum likelihood (ML), quasi likelihood (QL) and weighted least squares (WLS) estimators for proportional error nonlinear regression models. This work was triggered by an application in thermoluminescece (TL) sedimentary dating for which the liteature revealed another estimator similar to weighted least squares with the exception of observed responses used as weights. This estimator that we refer to as data weighted least squares (DWLS) is also included in the comparison.

We show that on the order $\sigma,$ all four estimators behave similar to
 ordinary least squares estimators for standard linear
regression models. On the order of $\sigma^2,$ the estimators have  biases. Formulae that are valid in the limit of small 
measurement error are derived for the biases and the variances of the four estimators. The maximum likelihood estimator has less bias compared to the quasi likelihood estimator. Conditions are derived under which weighted least squares and maximum likelihood estimators have similar biases.   
On the order of $\sigma^{2}$, all estimators have similar standard errors. On higher order of $\sigma$, the maximum likelihood estimator has smaller variance compared to the quasi likelihood estimator,  provided that the random errors have the same first four moments as the normal distribution.

The maximum likelihood and quasi-likelihood estimating equations are unbiased. In large samples, these two estimators are distributed as multivariate normal. The estimating equations for weighted least squares and data weighted least squares are biased. However,  in the limit of $\sigma  \to 0$ and $n \to \infty,$ if $ n^{1/2} \sigma$ remains bounded, these two estimators are also distributed as multivariate normal.
A simulation study justified the applicability of the derived formulae in the presence of measurement errors typical in sedimentary data.
Results  are illustrated with a data set from thermoluminescence sedimentary dating.  The work reported is applicable to more general contexts such as those arising in change point regression analysis.



\section{Introduction}
\label{introd}
Comparison of estimators based on large sample asymptotics is quite common in statistical literature. However, such comparisons are less appealing for contexts where sample sizes are relatively small. Kadane \cite{kadane} proposed comparison of estimators in the limit of small measurement errors  and reported that small sigma asymptotics can provide definite answers to normative choice of estimators. In this paper, we present some useful results from the comparison of estimators for nonlinear regression models with small measurement errors proportionately changing with the mean.

 This work was triggered by an application in thermoluminescence (TL) sedimentary dating in which typical data sets are small and have relatively small measurement errors. Apart from maximum likelihood, quasi likelihood and weighted least squares that are well known, literature on sedimentary data analysis reveals another estimator similar to weighted least squares with the exception of observed responses used as weights. We refer to this estimator as data weighted least squares (DWLS).  is also included for comparison. 

  In Section~\ref{models}, we present the notation and outline the estimating equations for  these estimators.   
In Section \ref{bias}, we derive formulae for the biases and 
variances of the estimators for theses models that are valid in the limit of small 
measurement error. We show that the maximum likelihood estimator has less bias compared to the quasi likelihood estimator. Conditions are derived under which weighted least squares and maximum likelihood estimators have similar biases.   
We further show that maximum likelihood estimators have smaller variances 
compared to quasi likelihood estimators, provided that the random errors 
have the same first four moments as the standard normal distribution.

Standard large sample small sigma distributional approximations 
for these estimators are presented in Section~\ref{dist}.  
The weighted and data weighted least squares estimators are not consistent in the 
limit of fixed measurement error. The trade-off between small
measurement error and bounds on the sample size needed to permit useful distributional
approximations are also examined in Section ~\ref{dist}. We show that the small $\sigma$  asymptotic results
remain relevant provided $\sigma\sqrt{n}$ is not large; the relative measurement error, $\sigma$ is scale free and the bounds on $\sigma$ are applicable in general.

In Section \ref{simulation}, we present the results of a Monte Carlo study that closely mimic an application in TL sedimentary data analysis.  The theoretical results derived in this 
article are demonstrated in Section \ref{example}, using a data set from sedimentary dating. 
Section \ref{conclude} offers some concluding remarks.

\section{Proportional error nonlinear regression model and estimating equations}
\label{models}
 
The proportional error nonlinear regression model that we focus in this study is
$y_{i}=f(x_{i},\vec{\theta}_{0}) (1 + \sigma \epsilon_{i}),$ 
where $\vec{\theta}_{0}$ denotes the vector of unknown true parameters, $\sigma$ denotes the relative error in a single measurement and the mean response $f(x_{i},\vec{\theta}_{0})$ is any nonlinear function. 
For notational convenience, we write $f_{i}$ for $f(x_{i},\vec{\theta})$ 
and omit the suffix $i$ when there is no confusion. Let $\hat{\vec{\theta}} $ denote an estimator for $\vec{\theta}$. 
Let $\bigtriangledown f(x_{i},\vec{\theta}) = 
      \partial f(x_{i},\vec{\theta})/\partial \vec{\theta}^{T}$ 
denote the gradient vector. 
Let $\bigtriangledown f_{o}$ and 
$\bigtriangledown f_{\hat{\vec{\theta}}} $ denote the gradient vector 
evaluated at $\vec{\theta}_{o}$ and $\hat{\vec{\theta}}$ respectively. Let $l$ denote the log-likelihood assuming normally distributed errors. Maximum likelihood estimator $\hat{\theta}$   solves the system of equations
\begin{eqnarray*}
 \left \{ \left. \frac{\partial l}{\partial \theta}\right 
   |_{\hat{\theta}, \hat{\sigma}}=0, \ \ 
  \left. \frac{\partial l}{\partial \sigma}\right 
   |_{\hat{\theta}, \hat{\sigma}}=0, \right \}
\end{eqnarray*}

The estimating equations for the four estimators are:

\begin{eqnarray}
\label{mle}
 & \text{ML:} & \frac{1}{n} \left\{ \sum_{i=1}^{n}
\frac{\left(y_{i}-{f}\right)^{2}}{\hat{f}^{2}} \right\} \sum_{i=1}^{n}{\frac{\bigtriangledown
\hat{f}
}{\hat{f}}  -  \sm \frac{(y_{i}-\hat{f}) }{{\hat{f}}^{2}} \bigtriangledown
\hat{f}-  \sm \frac{(y_{i}-\hat{f})^{2} }{{\hat{f}}^{3}} \bigtriangledown
\hat{f} } = 0. \\
& \text{QL:} &
\label{qle}
   \sum_{i=1}^{n}{ \frac{\left \{ y_{i} - f(x_{i},\hat{\theta})
\right \}
}{f_{\hat{\theta}}^{2}} \bigtriangledown f_{\hat{\theta}}}=0 \\
& \text{WLS:} &
\label{wls}
  \sum_{i=1}^{n}{\left \{ \frac{ \left ( y_{i}
- f_{
\hat{\theta}} \right ) }{f_{\hat{\theta}}^{2}} \bigtriangledown
f_{\hat{\theta}}
 \right \} } + \sum_{i=1}^{n}{\left \{ \frac{ \left ( y_{i} -
f_{\hat{\theta}} 
\right )^{2} }{f_{\hat{\theta}}^{3}} \bigtriangledown f_{\hat{\theta}}\right
\} }=0 \\
& \text{DWLS:} &
\label{dwls}
  \sum_{i=1}^{n}{ 
    \left[ \frac{\{y_{i} - f(x_{i},\hat{\theta})\}}{y_{i}^{2}} 
     \bigtriangledown f_{\hat{\theta}} \right] } =0
 \end{eqnarray}

We begin our analysis of small $\sigma$ asymptotics by 
approximating $\hat{\vec{\theta}}$ using the expansion 
$\vec{\that}=\vec{\theta}_{0}+C_{1} \sigma + C_{2} \sigma^{2},$ 
where $C_{1}$ and $C_{2}$ are $p \times 1$ random vectors 
that do not depend on $\sigma.$ 
Let  $H(x_{i},\vec{\theta}) = \partial^{2} f(x_{i},\vec{\theta})
        /\partial \vec{\theta} \partial \vec{\theta}^{T}$ 
denote the Hessian matrix. 
Further, let
$\bJ= \bigtriangledown f(x_{i},\vec{\theta}_o)/f(x_{i},\vec{\theta}_o)$
and $\bK = H(x_{i},\vec{\theta}_o)/f(x_{i},\vec{\theta}_o)$.  
For $\vec{\that}$ close to $\vec{\theta}_{0},$ the second order Taylor 
approximation for $f(x_{i},\hat{\vec{\theta}})$ around $\vec{\theta}_{0}$ 
can be written as:
$ f_{\hat{\vec{\theta}}}  \approx  
  f_{0} + (\hat{\vec{\theta}}-\vec{\theta}_{0})^{T} \bigtriangledown f_{0} 
  + \frac{1}{2} (\hat{\vec{\theta}}-\vec{\theta}_{0})^{T} H_{0} (\hat{\vec{\theta}}-\vec{\theta}_{0}),
$ where $H_{0}$ denote the Hessian matrix evaluated at $\vec{\theta}_{0}$.

Neglecting terms of $O(\sigma^{3})$ and higher, we find
\begin{equation*}
\label{fhat}
 f_{\hat{\vec{\theta}}} \approx  f_{0} 
  + (C_{1}^{T} \bigtriangledown f_{0}) \sigma 
  + (C_{2}^{T} \bigtriangledown f_{0} 
  + \frac{1}{2} C_{1}^{T} H_{0} C_{1}) \sigma^{2}.
\end{equation*}
Using this approximation, the estimating equations \ref{mle} to \ref{dwls} can be written as follows:

ML:
\begin{align*}
\label{Emle}
 & \frac{1}{n} \left\{ 
  \sum_{i=1}^{n} \left(\ei^{2} - 2 \ei C_{1}^{T} \bJ
 + C_{1}^{T} \bJ C_{1}^{T} \bJ \right) \sigma^{2} 
     \left(1-2C_{1}^{T}\bJ  \sigma\right) \right\} \\
& \qquad \times   \sum_{i=1}^{n}{ \left \{  \left (\bJ  
   + \bK C_{1} \sigma + \bK C_{2} \sigma^{2} \right ) 
     \left( 1 - C_{1}^{T} \bJ  \sigma 
     + C_{1}^{T} \bJ C_{1}^{T} \bJ  \sigma^{2} \right) \right \} } \\
 &  - \sum_{i=1}^{n}{ \left\{\left(\epsilon_{i} - C_{1}^{T} \bJ \right ) \sigma 
   - \left( C_{2}^{T} \bJ  + \frac{1}{2} C_{1}^{T} \bK C_{1} \right) \sigma^{2} 
     \right\} 
    \left ( 1 - 2 C_{1}^{T} \bJ  \sigma \right ) }
\left (\bJ + \bK C_{1} \sigma + \bK C_{2} \sigma^{2} \right )  \\
& - \sum_{i=1}^{n}{ \left \{ (\epsilon_{i}^{2} - 2 \epsilon_{i} C_{1}^{T} \bJ 
   + C_{1}^{T} \bJ C_{1}^{T} \bJ ) \sigma^{2} 
     \left( 1 - 3 C_{1}^{T} \bJ  \sigma \right) \right \} }
\left( \bJ + \bK C_{1} \sigma + \bK C_{2} \sigma^{2}\right) = 0. 
\end{align*}

QL:
\begin{eqnarray*}
\label{Eqle}
\sum_{i=1}^{n}{\left \{ ( \epsilon_{i} - C_{1}^{T}\bJ)
- (C_{2}^{T} \bJ
 + \frac{1}{2} C_{1}^{T} \bK C_{1} ) \sigma \right \} 
 \left \{ 1 - 2 C_{1}^{T}\bJ   \sigma \right \} } 
   (\bJ + \bK C_{1} \sigma + \bK C_{2} \sigma^{2} ) = 0
\end{eqnarray*}

WLS:

\begin{align*}
&\sum_{i=1}^{n}{ \left \{ (
\epsilon_{i} -
 C_{1}^{T} \bJ ) - (C_{2}^{T}
\bJ + \frac{1}{2} C_{1}^{T} \bK C_{1} ) \sigma \right\} 
 \left(1-2 C_{1}^{T} \bJ \sigma \right)}
 \left(\bJ + \bK C_{1} \sigma + \bK C_{2} \sigma^{2}\right) \\
& \quad  + \sum_{i=1}^{n} {\left( \epsilon_{i}^{2} - 2 \epsilon_{i}
C_{1}^{T} \bJ + C_{1}^{T}
\bJ  C_{1}^{T} \bJ \right) \sigma
\left( 1 - 3 C_{1}^{T} \bJ \sigma \right)} 
\left(\bJ  + \bK C_{1} \sigma + \bK C_{2} \sigma^{2}\right)=0.
\end{align*}

DWLS:

\begin{eqnarray*}
 \sum_{i=1}^{n} \left \{ 
  \left(\epsilon_{i} - C_{1}^{T} \bJ f_{0}\right) 
   - \left(C_{2}^{T} \bJ  + \frac{1}{2} C_{1}^{T} H_{0} C_{1} \right) 
     \sigma \right \} 
\left(1 - 2 \sigma \epsilon_{i} + 3 \sigma^{2} \epsilon_{i}^{2}\right)  
\left(\bJ + \bK C_{1} \sigma + \bK C_{2} \sigma^{2}\right) = 0.
\end{eqnarray*}

\subsection{Biases of the estimators}
\label{bias}

Several authors have discussed bias correction for special classes
of nonlinear regression models. For example, Box \cite{box} and Cook {\em et~al.} \cite{cook} have addressed the problem of computing the biases of the least squares
estimators for parameters in standard nonlinear regression models. 
Paula \cite{paula} has discussed bias correction to the order $O(1/n)$
for exponential family nonlinear models. We discuss bias correction in the estimators for nonlinear regression models in which standard deviation is proportional to the mean. Formulae are presented for the biases and standard errors  that are valid in the limit of small measurement error. Based on the formulae, some  useful asymptotic results for comparison of the biases of the four estimators are derived.
\begin{Res}
To order $\sigma,$ the estimators maximum likelihood, quasi likelihood, weighted least squares and data weighted least squares estimators 
 behave similar to ordinary least squares estimators in standard linear regression models.
\end{Res}

\noindent {\bf {\em Proof}}
Equating the coefficients of powers of 
$\sigma$, in the estimating equations to zero we find that
in all four estimation methods, $C_{1}$ can be written as
\[ C_{1}  =  \left( \sum_{i=1}^{n} \bJ \bJT \right)^{-1} 
  \left (\sum_{i=1}^{n} \epsilon_{i} \bJ  \right) =\left({\bf J}^T{\bf J}\right)^{-1} {\bf J}^{T}\vec{\epsilon},  \]
where $\bf J$ is the $n \times p$ matrix with $J_i^T$  as the $i$th
row and $\vec{\epsilon}$ is the $n \times 1$ vector with entries
$\epsilon_i$. 

Thus, to order $\sigma$ the standardized
estimation error $\frac{\vec{\hat{\theta}}-\vec{\theta}_{0}}{\sigma} $ has the form
$$
\frac{\that- \theta_o}{\sigma} = \left({\bf J}^T{\bf J}\right)^{-1}
{\bf J}^T \boldsymbol\epsilon
$$
 This is the usual ordinary least squares
formula for a regression problem with design matrix $\bf J$.
and hence the result  follows for general nonlinear regression models with proportional errors.

\subsection{Standard errors of the estimators}
\label{osigma2}
Formulae for the biases and the standard errors of the four estimators on the order $O(\sigma^{2})$ can be derived by considering coefficients of $\sigma^{2}$ in the four estimating equations presented in Section ~\ref{models}. With some algebra (see ~\cite{perera} for details), the term $C_{2}$  in the estimating equations can be written in the form 
$\left({\bf J}^T{\bf J}\right)^{-1} A$,  
where the random error term $A$  for each method is presented in Table ~\ref{Aterm}. 
\begin{table}[h]
\begin{center}
\begin{tabular}{|l|l|} \hline
Method of &  \hspace{3cm} $A$\\
Estimation  &  \\ \hline
ML & $ 
  -  \frac{1}{n} \left [ \sum_{i=1}^{n} \left \{\ei^{2}-2 \ei C_{1}^{T} \bJ  
     + C_{1}^{T} \bJ C_{1}^{T} \bJ  \right \} \right ] 
      \left \{ \sum_{i=1}^{n}{\bJ  } \right \} 
       +  \sum_{i=1}^{n}{ \bK C_{1} \epsilon_{i}} 
 - 4 \sum_{i=1}^{n}{ \bJ  \bJT C_{1} \epsilon_{i} }$  \\
 & $ \quad - \sum_{i=1}^{n}{C_{1}^{T} \bJ  \bK C_{1}} 
    +  3 \sum_{i=1}^{n}{C_{1}^{T} \bJ   \bJT C_{1} \bJ  }
 - \frac{1}{2} \sum_{i=1}^{n}{C_{1}^{T} \bK C_{1} \bJ } 
  + \sum_{i=1}^{n}{\epsilon_{i}^{2} \bJ  }.$  \\ \hline
QL & $\sum_{i=1}^{n}{\bK \epsilon_{i} C_{1} } 
   - 2   \sum_{i=1}^{n}{\bJ \bJT C_{1} \epsilon_{i} } 
   - \sum_{i=1}^{n}{C_{1}^{T} \bJ \bK  C_{1} }$ \\
 & $ + 2 \sum_{i=1}^{n}{C_{1}^{T} \bJ \bJT C_{1} \bJ } 
   - \frac{1}{2} \sum_{i=1}^{n}{C_{1}^{T}\bK C_{1} \bJ }$ \\ \hline
WLS & $\sum_{i=1}^{n}{\epsilon_{i}\bK  C_{1} } 
    - 4 \sum_{i=1}^{n}{ \bJ \bJT C_{1} \epsilon_{i}} 
    - \sum_{i=1}^{n}{C_{1}^{T} \bJ\bK C_{1}} 
+ 3 \sum_{i=1}^{n} C_{1}^{T} \bJ C_{1}^{T} \bJ \bJ$  \\
  &  $- \frac{1}{2} \sum_{i=1}^{n} C_{1}^{T} \bK  C_{1} \bJ 
  + \sum_{i=1}^{n} \epsilon_{i}^{2}  \bJ $  \\ \hline
DWLS & $\sm \ei \bK C_{1} 
  - 2 \sm \ei^{2} \bJ- \sm C_{1}^{T} \bJ \bK C_{1} 
 + 2 \sm \ei C_{1}^{T} \bJ\bJ
   - \frac{1}{2} \sm C_{1}^{T} \bK C_{1} \bJ $   \\ \hline
\end{tabular}
\end{center}
\caption{The random error term $A$ contributing to the bias on the order $\sigma^{2}$}
\label{Aterm}
\end{table}

With some algebra, we derived the formulae  presented in Table~\ref{bias.est} for the biases and variances,  where we use the notation
\begin{eqnarray*}
p & = & \mbox{Number of components of $\theta$} \\
w_{1,i} &=& tr \left \{ \bJ  \bJT \left({\bf J}^T{\bf J}\right)^{-1}
\right \} \\
\mbox{and} \ \ w_{2,i} & = & tr \left \{ \bK \left({\bf J}^T{\bf
J}\right)^{-1}
\right \}.
\end{eqnarray*}
Notice that the matrix $J$ plays the role of the design matrix in the general linear regression model and $w_{1,i}$ are the diagonal entries in the corresponding ``hat'' matrix
${\bf J} \left({\bf J}^T{\bf J}\right)^{-1}{\bf J}^T$. 
These formulae permit us to make the useful observation that, on the order $O(\sigma^{2}),$ all four estimators have the same standard error. Therefore, biases on the order $O(\sigma^{2})$ are useful in choosing between these four estimators. The fact that the weights have to be updated at each iteration makes obtaining data weighted least squares estimates computationally much simpler compared to weighted least squares.
\begin{table}[htbp]
\begin{center}
\begin{tabular}{|l|l|r|} \hline
Method of &  \hspace{3cm} Bias & $Var(\that)$ \\
Estimation  & & \\ \hline
ML & $\left({\bf J}^T{\bf J}\right)^{-1} \left \{ - \sm (w_{1,i}-\frac{p}{n}) \bJ
   - \frac{1}{2} \sm w_{2,i} \bJ\right \} \sigma^{2}$ 
    & $\sigma^{2}\left({\bf J}^T{\bf J}\right)^{-1}$ 
\\ & & \\ \hline
QL & $\left({\bf J}^T{\bf J}\right)^{-1}
\left \{ - \frac{1}{2} \sm w_{2,i} \bJ\right \} \sigma^{2}$ & $
 \sigma^{2}\left({\bf J}^T{\bf J}\right)^{-1}$ 
\\ & &  \\ \hline
WLS & $\left({\bf J}^T{\bf J}\right)^{-1}
\left \{ \sm \bJ- \sm w_{1,i} \bJ- \frac{1}{2} \sm w_{2,i} \bJ
\right \} \sigma^{2}$ & $ \sigma^{2}\left({\bf J}^T{\bf J}\right)^{-1}$
\\ & &  \\ \hline
DWLS & $\left({\bf J}^T{\bf J}\right)^{-1}
\left \{ -2 \sm \bJ+ 2 \sm w_{1,i} \bJ- \frac{1}{2} \sm w_{2,
i} \bJ\right \} \sigma^{2}$ & $ \sigma^{2}\left({\bf J}^T{\bf J}\right)^{-1}$
\\ & &  \\ \hline
\end{tabular}
\end{center}
\caption{The biases and the variances of the estimators}
\label{bias.est}
\end{table}

\begin{Res}
In normal error nonlinear regression models of the form $y= f(x,\theta_{1},\cdots,\theta_{p})(1+\sigma \epsilon),$ if the response function can be written as
$f(x,\theta_{1},\cdots,\theta_{p})=\theta_{1}
 f^{*}(\theta_{2},\cdots,\theta_{p})$, where $f^{*}$ is some function that
does not depend on $\theta_{1},$  in the limit of small measurement errors,  maximum likelihood estimators  and weighted least squares estimators for all the parameters except $\theta_{1}$ have identical biases. 
\label{res1}
\end{Res}

\noindent {\em {\bf Proof}}\/: 
Let  $\bigtriangledown f^{*}$ be the gradient vector of length $p-1$ consisting of the derivatives with
respect to $\theta_{2},\cdots,\theta_{p}.$  Let  $J_{i}^{*} = \bigtriangledown f^{*}/f^{*}$  .
For 
response functions of the form
considered here,
it is easy to see that
$  \bJ $ is of the form
\[  \bJ 
 = \left( \begin{array}{ll}
                       \frac{1}{\theta_{1}} \\
                        \bJstar 
                       \end{array} \right ). \]
and
\[ {\bf J}^T{\bf J}= \left [ \sum_{i=1}^{n}{ \bJ \bJT
 }\right ] 
   = \left[ \begin{array}{lll}
                      \frac{n}{\theta_{1}^{2}} &  \frac{1}{\theta_{1}} \sm
 \bJstarT \\
\frac{1}{\theta_{1}} \sm \bJstar
& \sm  \bJstar \bJstarT
\end{array}
\right]. \]

Thus,
$\left({\bf J}^T{\bf J}\right)^{-1}\sm \bJ $ 
takes the form $ [\theta_{1}, 0, \ldots, 0]^{T}.$ The result  immediately follows from
 the formulae presented in Table~\ref{bias.est}.

\section{Large sample small sigma behaviour of the estimators}
\label{dist}

In this section, we provide large sample distributional approximations 
 that are valid for general proportional error nonlinear regression models in the limit of small measurement errors. First note that all four estimators are defined as roots of a
general estimating equation of the form
\begin{eqnarray*}
H_{n}(\theta) = \frac{1}{\sqrt{n}} \sm h_{i}(y_{i},\theta) = 0,
\end{eqnarray*}
where $h_{i}(y_{i},\theta)$ is a function of $y_{i}$ and $\theta$.
We can study the large sample behaviour as usual by studying $H_n$. In what
follows, we are assuming standard regularity conditions such as
$
\sum \bJ = O(n)$, $\sum \bJ\bJT=O(n)$ and similar conditions on the second
derivatives.

\subsection{Large sample small $\sigma$ behaviour of ML}
If ${ E}(Y_{i})=f(x_{i},\vec{\theta})$ and
${Var}(Y_{i})=\sigma^{2}f^{2}(x_{i},\vec{\theta})$ and the assumed error distribution is correct, the maximum likelihood estimating equations are unbiased.
In large samples, assuming ${ E}(Y_{i}^{4}) < \infty$ we find
\[ \left( \begin{array}{lll}
         \that - \theta \\
         \hat{\sigma} - \sigma
         \end{array}
         \right) \sim MVN \left (0, \;
{E}\left[-H_{n}^{\prime}(\theta)\right]^{-1
} {\rm Var} \left\{H_{n}(\theta)\right\}{ E}\left[-H_{n}^{\prime}(\vec{\theta})\right]^{-1}
\right ), \]
where 
$$
 H_{n}^{\prime}(\vec{\theta}) = \left [ \begin{array}{lll}
                           \frac{\partial^{2} l}{\partial \vec{\theta} \partial
\vec{\theta}
^{T}} & \frac{\partial^{2} l}{\partial \theta \partial \sigma} \\
\frac{\partial^{2} l}{\partial \vec{\theta}^{T} \partial \sigma} &
\frac{\partial^{2}
l}{\partial \sigma^{2}}
\end{array}
\right ] 
$$
and $l$ denotes the log-likelihood. 
We show that the variance covariance matrix of the maximum 
likelihood estimator $\hat{\theta}$ (see Appendix) is:
\[ V(\that) = \left[
              \left( 2+\sigma^{-2} \right) \sm \bJ\bJ^T
- 2n^{-1}\left( \sm \bJ \right) \left( \sm \bJ \right)^T
\right
]^{-1}. \]
Writing $n^{-1} \sm \bJ = \overline{J}$ the above covariance matrix can be rewritten as (see Appendix)
\begin{equation}
\label{var.that.dist}
V(\that) =\sigma^2 \left[
              {\bf J}^T{\bf J} + 2 \sigma^2\left\{ \sm
\left(
 \bJ-\overline{J}\right)\left(\bJ-\overline{J}\right)^{T}
\right\}
\right]^{-1}.
\end{equation}

\subsection{Large sample small $\sigma$ behaviour of QL}
In the case of quasi likelihood,
$$
H_{n}(\theta)=n^{-1/2} \sm \left \{ \frac{y_{i}-f(x_{i},\theta)}{f_{i}^{2}}
\right \} \trgfi. 
$$
Clearly, under the assumptions that $E(y_i)=f(x_i,\theta)$ 
(as is implied by our model) and ${\rm Var}(y_i)< \infty$ we have
${\rm E}\{h_{i}(y_{i},\theta)\}=0$ and ${\rm Var}\{h_{i}(y_{i},\theta)\}<\infty $ .
Hence, quasi-likelihood estimating equations are unbiased.
General large sample considerations will then establish that, in large
samples,
\begin{align*} (\that -\theta) 
& \approx \left\{-H_n^\prime(\theta)\right\}^{-1}H_n(\theta) +o_p(n^{-1/2})
\\
&  \sim MVN
\left(0,E\left[-H_{n}^{\prime}(\theta)\right]^{-
1} Var \left(H_{n}(\theta)\right)
E\left[-H_{n}^{\prime}(\theta)\right]^{-1} 
\right ), 
\end{align*}
where 
$H_{n}^{\prime}(\theta)$ is the derivative of $H_{n}(\theta)$ with
respect to $\theta$.  
It is easy to see that
$ E\left[H_{n}^{\prime}(\theta)\right]=\sm \bJ \bJ^{T}$ so that
$E\left[H_{n}^{\prime}(\theta)\right]^{-1}= \left[\sm \bJ
\bJ^{T}\right]^{-1} = \left({\bf J}^T{\bf J}\right)^{-1}$ and
\[ {Var}\left(H_{n}(\theta)\right)=\sm \frac{Var(Y_{i})}{f_{i}^{2}} \bJ\bJ^{T}.\]
Thus, in large samples, 
\[ (\that -\theta) \sim MVN \left(0,\left({\bf J}^T{\bf J}\right)^{-1}\left[ \sum_{i=1}^{n}
\frac{Var(Y
_{i})}{f_{i}^{2}} \bJ \bJ^{T} \right] \left({\bf J}^T{\bf J}\right)^{-1} \right). \]
For models with ${\rm Var}(Y_{i}) = \sigma^{2} f^{2}(x_{i},\theta),$  
the asymptotic expansion above simplifies to give
the result that
\[ (\that -\theta) \sim MVN \left(0, \sigma^{2} \left({\bf J}^T{\bf J}\right)^{-1}\right). \]

\begin{Res}
In proportional error nonlinear regression models with normally distributed errors, in the limit of small $\sigma,$ the maximum likelihood estimator has smaller variance compared to quasi-likelihood estimators.
\end{Res}

\noindent {\em {\bf Proof}}\/: 
The term $ \sm \left(\bJ-\overline{J}\right)\left( \bJ-\overline{J}\right)^{T}$ is positive definite. Thus, the result follows immediately from the formulae derived for the variances of the two estimators.
 
\begin{Res}
In large samples with small measurement errors, on the order $o(\sigma^{4}),$ the maximum likelihood estimator for proportional error nonlinear regression models have the same variance as in the estimators for normal error general linear regression models.
\end{Res}

\noindent {\em {\bf Proof}}\/:  
 On the order $o(\sigma),$  ignoring the second term in Equation
\ref{var.that.dist},
  we find that the  variance covariance matrix reduces to the
variance
 covariance matrix for the general linear regression model with the design matrix replaced with the matrix $J.$

\subsection{Large sample small $\sigma$ behaviour of WLS and DWLS}

The estimating equations for weighted and data weighted least squares are biased. Thus, these estimates are not consistent as $n\to\infty$ with $\sigma$ fixed.  If we consider a limit in which $\sigma\to 0$
as $n\to\infty$ in such a way that $n^{1/2}\sigma$ remains bounded, then we may get normal approximations.  To simplify the presentation we assume that the following limits exist:
\begin{align*}
n^{1/2}\sigma & \to \delta
\\
({\bf J}^T{\bf J}/n)^{-1}& \to \boldsymbol\Sigma
\\
\sum \bJ/n & \to \boldsymbol\Gamma_1
\\
\sum w_{1,i} \bJ/n &  \to \boldsymbol\Gamma_2
\\
\sum w_{2,i} \bJ/n &  \to \boldsymbol\Gamma_3
\end{align*}
Under these conditions, we find that the limiting distribution for weighted least squares is
$$
\frac{\sqrt{n}(\that-\theta)}{\sigma} \Rightarrow MVN\left(\delta \boldsymbol\Sigma
 (\boldsymbol\Gamma_1-\boldsymbol\Gamma_2 - \boldsymbol\Gamma_3/2),\boldsymbol\Sigma
 \right)
 $$
 and that for data weighted least squares is
$$
\frac{\sqrt{n}(\that-\theta)}{\sigma} \Rightarrow MVN\left(\delta \boldsymbol\Sigma
 (-2 \boldsymbol\Gamma_1+2\boldsymbol\Gamma_2 - \boldsymbol\Gamma_3/2),\boldsymbol\Sigma
 \right).
 $$
 A further level of approximation can be noted.  The weights $w_{1,i}$ and $w_{2,i}$
 have a sum over $i$ which should be $O(1)$.  This means that usually we will have
 $$
 \boldsymbol\Gamma_2 = \boldsymbol\Gamma_3 = 0.
 $$
Our distributional approximations then simplify to give the following asymptotic results for WLS and DWLS, in the limit of large $n$ and small $\sigma$ such that $n^{1/2} \sigma$ is bounded:

\begin{center}
WLS:  \hspace{1cm} $\frac{\sqrt{n}(\vec{\that}-\vec{\theta})}{\sigma} \Rightarrow MVN\left(\delta \boldsymbol\Sigma
 \boldsymbol\Gamma_1,\boldsymbol\Sigma  \right)  $ \\
 DWLS:  \hspace{1cm} $\frac{\sqrt{n}(\vec{\that}-\vec{\theta})}{\sigma} \Rightarrow MVN\left(-2\delta \boldsymbol\Sigma  \boldsymbol\Gamma_1,\boldsymbol\Sigma  \right)  $  \\
\end{center}

In passing we also note that for mean functions $f(x,\vec{\theta})$ such as the saturating exponential model, according to Result \ref{res1}, $ \boldsymbol\Gamma_{1}$ will have all but the first entry 0.

\section{Simulation study}
\label{simulation}
Now we describe the results of a simulaiton study that examine the finite sample applicability of the derived asymptotic results. The simulation study mimic an application in TL sedimentary dating using an experimental design called the partial bleach method. More simulation results based on other experimental designs used in TL studies are presented in Perera \cite{perera}.  In the partial bleach method, the sediments are dated based on an estimate for what is known as the equivalent dose. Mathematically, the equivalent dose is the absolute value of the dose level, $x$ corresponding to the point of intersection of two nonlinear functions fitted for two data sets known as unbleached data and bleached data. The functions fitted are:
$f_{1}(x,\vec{\theta}_{1})= \alpha_{1}\left( 1- \exp \left( -\frac{x+\alpha_{2}}{\alpha_{3}}\right)\right)$,  where $\vec{\theta}_{1}= (\alpha_{1},\alpha_{2},\alpha_{3})^{T}$ 
and
$f_{2}(x,\vec{\theta}_{2})= \beta_{1}\left( 1- \exp \left( -\frac{x+\beta_{2}}{\beta_{3}}\right)\right)$,  where $\vec{\theta}_{2}= (\beta_{1},\beta_{2},\beta_{3})$. Let $\vec{\theta}=(\vec{\theta}_{1},\vec{\theta}_{2})^{T}.$ The equivalent dose $\gamma$ is estimated as a root of the equation $g(x,\vec{\theta})=f_{1}(x,\vec{\theta}_{1})-f_{2}(x,\vec{\theta}_{2})=0.$

 For the simulation study, dose levels and sample sizes were fixed in advance at the levels in QNL84-2  experimental data set proposed by  Berger {\em et~al.} \cite{hun1}. The sample sizes  of unbleached and bleached data sets were $n_{1}=16$ and $n_{2}=13$ respectively. Compared to the number of fitted parameters, sample sizes are relatively small.  The TL intensity $y$ was generated according to $y=f(x,\vec{\theta}_{j})(1+\sigma \epsilon),$ for $j=1,2$ by setting the parameter values at the maximum likelihood estimates obtained for the QNL84-2 data. Thus, we assigned
$\alpha_{1}=142853.0, \alpha_{2}=123.182, \alpha_{3}=393.065,\beta_{2}=192.547 $ and $\beta_{3}=756.620$.  The parameter $\gamma$ corresponding to the equivalent dose was set at $\gamma=-87.45$; since curves intersect over the region of negative $x$, this correspond to an equivalent dose of 87.45 Gray. The value of $\beta_{1}$ was taken to be $\beta_{1}=\frac{\alpha_{1} \left( 1- \exp \left( -\frac{\gamma+\alpha_{2}}{\alpha_{3}}\right)\right)}{\left( 1- \exp \left( -\frac{\gamma+\beta_{2}}{\beta_{3}}\right)\right)}$ so that the two curves are guaranteed to intersect at $\gamma$. The values of $\sigma$ chosen common to both curves, biases computed using the derived formulae ($B_{T}$) and the estimated biases based on 10000 simulations for each case ($B_{s}$) are presented in Table \ref{results}. 

\begin{table}[h]
\begin{center}
\begin{tabular}{|c|c|c|c|c|c|c|c|c|} \hline
$\sigma$ & \multicolumn{2}{c|}{ML} & \multicolumn{2}{c|} {QL}& \multicolumn{2}{c|} {WLS} & \multicolumn{2}{c|} {DWLS} \\
  \cline{2-9} 
 & $B_{T}$  & $B_{s}$  & $B_{T}$  & $B_{s}$ & $B_{T}$  & $B_{s}$ & $B_{T}$  & $B_{s}$ \\  \hline
0.01 &- 0.046 & -0.046 &-0.049 & -0.048 & -0.046 & -0.045 & -0.054 & -0.045 \\
 0.02 & -0.182 & -0.181 & -0.195 & -0.195 & -0.183 & -0.182 & -0.217 & -0.221 \\ 
0.03 & -0.410 & -0.429 & -0.438 & -0.444& -0.412 & -0.414 & -0.489 & -0.508 \\ 
0.04 & -0.730 & -0.783 & -0.778 & -0.824 & -0.733 & -0.784 & -0.869 & -0.923 \\ 
0.05 & -1.140 & -1.289 & -1.216 & -1.329 & -1.146 & -1.267 & -1.358 & -1.483 \\ 
0.06 & -1.641 & -1.779  & -1.752 & -1.865 & -1.650 & -1.760 & -1.955 & -2.687 \\ \hline
\end{tabular}
\end{center}
\caption{Comparison of biases using the formulae and from simulation}
\label{results}
\end{table}

 The results  indicate good agreement between the biases computed from the derived formulae with the relevant biases estimated from the simulation study. 
We emphasize that as noted in Result \ref{res1}, both maximum likelihood and weighted least squares estimators for the equivalent dose have similar bias.

\section{Worked example}
\label{example}
We now use the QNL84-2 data set for further illustration of the derived results. 
The models described in Section \ref{simulation} were fitted to the data assuming a common $\sigma$. We note that fitting different $\sigma$ values for the two data sets gave similar parameter estimates (see Perera \cite{perera}).The biases and the mean squared errors (MSE) in Table \ref{exam} were estimated using the formulae given in Table \ref{bias.est} with parameters replaced by the corresponding estimates.  For maximum likelihood, we have used the maximum 
likelihood estimate for $\sigma$.  For the other three  methods we have used the 
unbiased estimate for $\sigma$ from the relevant fits.

\begin{table}[bhtp]
\begin{center}
\begin{tabular}{|c|c|c|c|c|c|c|} \hline
Data         & para. & Description  & \multicolumn{4}{c|}{Method} \\
\cline{4-
7}
             &       &  & ML & QL & WLS & DWLS \\  \hline

 QNL84-2       & $\alpha_{1} \times 10^{-4}$ & Estimate  &14.28 & 14.28 &14.30
& 14.25 \\
($n_{1}=16$) &              & bias & 0.02 & 0.03 & 0.05 & 0.09
\\
($n_{2}=13$)       &    & std. error (se) & 0.49 & 0.55 & 0.55 &
0.55 \\
      & & $bias/\sqrt{MSE} \times 100 \% $&4.08 &5.45&9.05 & 16.15 \\
\cline{2-
7}
             & $\alpha_{2}$ &  Estimate & 123.18 & 122.74 & 123.18 &
121.86 \\
             &        & bias     & 0.12    & 0.24   & 0.15   & 0.41  \\
             &  & std.error (se) &7.26& 8.12   & 8.16   & 8.10  \\
& & $bias/\sqrt{MSE} \times 100 \% $  &1.65& 2.95   & 1.84 & 5.06 \\
\cline{
2-7}
             & $\alpha_{3}$ & Estimate & 393.07  & 392.00 & 393.07
&389.92 \\
                &         & bias     & 1.64  & 2.46 & 2.07 &3.23   \\
                &   & std.error (se) & 33.11 & 37.04 & 37.20  &36.94  \\
& & $bias/\sqrt{MSE} \times 100 \% $  & 4.95  & 6.63 &5.56 &8.71 \\
\cline{2-7
}
         & $\beta_{2}$ & Estimate & 192.55  &193.37  & 192.54 & 195.18
\\
             &       & bias         & 0.39  & 0.72 & 0.49 & 1.19 \\
             & & std.error (se)     & 13.97 & 15.80&15.69 &16.12 \\
&& $bias/\sqrt{MSE} \times 100 \% $  &2.79 & 4.55 & 3.12 & 7.36  \\
\cline{2-7}
             & $\beta_{3}$ & Estimate & 756.62  & 761.65 & 756.59 &772.76
\\
             &     & bias           & 11.20 & 16.21 & 14.12  & 20.63  \\
           & & std.error (se)       & 105.46& 120.06 &118.49 & 124.19 \\
&&$bias/\sqrt{MSE} \times 100 \% $  & 10.56& 13.38 &11.83  & 16.39 \\ \cline{2-7}
& $\gamma$ & Estimate & 87.15  & 86.43 & 87.16 & 84.98
\\
             &     & bias           & 0.55 & 0.72 & 0.70 & 0.77  \\
           & & std.error (se)       & 9.13 & 10.14 & 10.26 & 9.97  \\
&&$bias/\sqrt{MSE} \times 100 \% $  & 6.01 & 7.08 & 6.81 & 7.70 \\ \cline{2-7}
   &  $\sigma$    & Estimate & 0.039   & 0.035 & 0.035 & 0.035  \\
\hline
\end{tabular}
\end{center}
\caption{Parameter estimates for the   QNL84-2  data set}
\label{exam}
\end{table}

The results of the worked example exemplifies that $\sigma$ is small as typical for sedimentary data. Furthermore, for all parameter estimates,  the relative biases are 
small compared to the standard errors.

\section{Concluding remarks and Discussion}
\label{conclude}
In this article, we focused on small relative measurement error 
asymptotics for maximum likelihood, quasi likelihood, weighted least 
squares and data weighted least squares estimators for parameters in  nonlinear regression models. Formulae valid in the limit of small 
measurement error were provided for the biases and mean squared errors 
of these estimators. Biases of maximum likelihood estimators 
were found to be smaller than the biases of quasi likelihood estimators. However, for 
certain parameters in specific models (see Result ~\ref{res1}), the biases of 
weighted least squares estimators were found to be similar to the 
biases of maximum likelihood estimators. Large sample asymptotics were presented for the four estimators and finite sample performance in the estimators were examined using simulations. The work was illustrated using the experimental data presented in Berger {\em et. al.} ~\cite{berg1}.

 The  work reported here has wider applications especially in the context of change point regression analysis. In contexts such as change point regression analysis, often one has to decide on whether a common relative error parameter $\sigma$  or different relative error parameters need to be fitted for  different segments. Intuitively, one should expect the biases and the standard errors of the estimators to depend on this decision. The DWLS estimating equations (see Section ~\ref{models}) for the proportionate error nonlinear models do not involve $\sigma$. Therefore, DWLS estimates are unchanged regardless of whether a common $\sigma$ or different $\sigma$'s are fitted for different segments. The estimating equations for the other three methods involve $\sigma$. For instance, quasi likelihood estimating equations for simultaneous curve fitting of two curves $f_{1}$ and $f_{2}$  for two segments with different $\sigma$'s take the form
$   \sum_{i=1}^{n_{1}}{ \frac{\left \{ y_{i} - f_{1}(x_{i},\hat{\theta})\right \}}{\sigma_{1}^{2} f_{1}(x,{\hat{\theta}})^{2}} \bigtriangledown f_{1}(x,{\hat{\theta})}}+
 \sum_{i=1}^{n_{2}}{ \frac{\left \{ y_{i} - f_{2}(x_{i},\hat{\theta})\right \}}{\sigma_{2}^{2} f_{1}(x,{\hat{\theta}})^{2}} \bigtriangledown f_{2}(x,{\hat{\theta})}}=0.$
Two-part iterative algorithms, each time solving estimating equations for $\vec{\theta}$ and  upgrading $\sigma_{1}$ and $\sigma_{2}$ using current parameter estimates need to be employed to estimate $\vec{\theta}$. If a common $\sigma$ is to be fitted, the curves have to be fitted simultanesously. Therefore, it is intuitive to expect that the parameter estimates for $\theta$ to depend on how we estimate $\sigma$'s. However, in contrary to what one expects,  for proportional error  nonlinear regression models, for fixed $\sigma$'s, the QL and WLS estimating equations for $\vec{\theta}$ are derivatives of a function (the likelihood for the gamma model or weighted error sum of squares) which is being optimized. The location of the optimum is invariant under reparametrization of $\vec{\theta}.$ When the curves are fitted separately for different segments, the estimating equations for QL and WLS clearly do not involve $\sigma$. Therefore, the invariance propoerty guarantees that the estimates for $\theta$ not to depend on whether we estimate $\sigma$  using the maximum likelihood estimate or using the least squares estimates.

Turning to ML, the situation is different. As for QL and WLS, the invariance property guarantees that the simultaneous curve fitting and separate curve fitting to yield same estimates. However, since maximum likelihood estimating equations are coupled with the estimating equations for $\sigma$'s, the estimates for $\vec{\theta}$ depend on how we estimate $\sigma$ and on whether the curves are fitted simultaneously or separately.

\section*{Appendix}
Here we prove that if the responses $Y_{i}$ have mean 
$f(x_{i},\theta)$, variance $\sigma^{2} f(x_{i},\theta)$,  and
\begin{align*}
{\rm E}\left\{\left(Y_{i}-f_{i}\right)^{3} \right\}& =0 \\
{\rm E}\left\{ \left(Y_{i}-f_{i}\right)^{4} \right\} &= 3 \sigma^{4} f_{i}^{4}
\end{align*}
then 
in large samples, the variance of the maximum likelihood 
estimator $\hat{\theta}$ is given by~(\ref{var.that.dist}).  Notice that if the errors in our
model have normal distributions then these assumptions on the third and fourth moments hold.

\noindent Proof:
In Section \ref{dist}, we noted that the variance of the maximum 
likelihood estimator is given by
$$
\left ( E\left[-H_{n}^{\prime}(\theta)\right]^{-1} 
  {\rm Var} \left(H_{n}(\theta)\right) E\left[-H_{n}^{\prime}(\theta)\right]^{-1}
\right ), 
$$
where
$$ H_{n}^{\prime}(\theta) = \left [ \begin{array}{lll}
                           \frac{\partial^{2} l}{\partial \theta \partial
\theta
^{T}} & \frac{\partial^{2} l}{\partial \theta \partial \sigma} \\
\frac{\partial^{2} l}{\partial \theta^{T} \partial \sigma} &
\frac{\partial^{2}
l}{\partial \sigma^{2}}
\end{array}
\right ] 
$$ 
and $l$ denotes the log-likelihood. 

Differentiating the log-likelihood function we find that, 
$E\left (- \frac{\partial^{2}
l}{\partial 
\theta \partial \theta^{T}}\right)$ can be written as $D^{T}MD$ 
where $D$ is the $n \times p$ matrix 
with $(i,j)$th entry $\partial f_{i} / \partial \theta_{j}$ and $M$ is
the diagonal matrix with $i$th diagonal element
$2/f_{i}^{2}+1/(\sigma^{2}f_{i}^{2})
$. So, $D^{T}MD$ can be written as
$\left(2+ \sigma^{-2}\right)\left( \sm \bJ
\bJ^{T}
\right)$. 
Now consider
\[ \frac{\partial^{2} l}{\partial \theta \partial \sigma} = -
\frac{2}{\sigma^{3
}} \sm \left(\frac{y_{i}-f_{i}}{f_{i}} \right) \bJ -
\frac{2}{\sigma^{3}} \sm
 \left(\frac{y_{i}-f_{i}}{f_{i}} \right )^{2} \bJ. \]
It is easy to see that \[ \; E \left[ - \frac{\partial^{2} l}{\partial
\theta 
\partial \sigma}\right]= \frac{2}{\sigma}\sm \bJ.\] \\
Since $ \frac{\partial^{2} l}{\partial \sigma^{2}} = \frac{n}{\sigma^{2}} -
\frac{3}{\sigma^{4}} \sm \left(\frac{y_{i}-f_{i}}{f_{i}} \right)^{2}, \;$ we
find $
\; E \left( - \frac{\partial^{2} l}{\partial \sigma^{2}}\right) = \frac{2
n}
{\sigma^{2}}.$ 
Now consider
\[ {\rm Var}\left\{H_{n}(\theta)\right\}= \left [ \begin{array}{lll}
                       {\rm Var} \left(\frac{\partial l}{\partial \theta}\right)
&
           {\rm  Cov}\left( \frac{\partial l}{\partial \theta},\frac{\partial
l}
{\partial \sigma}\right) \\
 \left\{{\rm Cov}\left( \frac{\partial l}{\partial \theta},\frac{\partial
l}{\partial 
\sigma}\right)\right\}^{T} &{\rm Var} \left(\frac{\partial l}{\partial
\sigma}\right)
 \end{array}
    \right ]. \]
The components of $Var \left(H_{n}(\theta) \right)$ can be computed as
follows:
\begin{align*}
{\rm Var} \left(\frac{\partial l}{\partial \theta}\right) & =
\sigma^{-2} \sm \bJ\bJ^{T} + \sigma^{-4} \sm
{\rm Var} 
\left( \frac{y_{i}-f_{i}}{f_{i}}\right)^{2} \bJ\bJ^{T}  + 
2 \sigma^{-4} \sm E \left[ \left(
\frac{y_{i}-f_{i}}{f_{i}}\right)^
{3} \right] \bJ\bJ^{T}, \\
\\
{\rm Cov}\left( \frac{\partial l}{\partial \theta},\frac{\partial l}{\partial
\sigma}
\right) & = \sigma^{-5}\sm E \left[  \left(
\frac{y_{i}-f_{i}}{f_{i}}
\right)^{3} \right ] \bJ +  \frac{1}{\sigma^{5}} \sm E \left[  \left(
\frac{y
_{i}-f_{i}}{f_{i}}\right)^{4} \right ] \bJ   - \frac{1}{\sigma} \sm \bJ,
\end{align*}
and
\begin{eqnarray*}
 Var \left(\frac{\partial l}{\partial \sigma}\right) =
\frac{1}{\sigma^{6}} \sm
 Var \left [ \left( \frac{y_{i}-f_{i}}{f_{i}}\right)^{2} \right].
\end{eqnarray*}
Now using our assumptions about the
third
 and fourth moments the components of ${\rm Var}\left\{H_{n}(\theta)\right\}$
simplify
 to give
\begin{align*}
{\rm Var} \left(\frac{\partial l}{\partial \theta}\right) & = 
\frac{1}{\sigma^{2}} 
\sm \bJ \bJ^{T} + 2 \sm \bJ \bJ^{T} \\
& =  \left( 2+\sigma^{-2}\right) \sm \bJ
\bJ^{T}, \\
{\rm Cov}\left( \frac{\partial l}{\partial \theta},\frac{\partial l}{\partial
\sigma}
\right) & =  \frac{2}{\sigma} \sm \bJ, \\
\mbox{and} 
 \ \ \ \ \ {\rm Var} \left(\frac{\partial l}{\partial \sigma}\right) & =  \frac{2
n}{\sigma^{
2}}.
\end{align*}
Thus we find that the usual Bartlett identity, ${\rm E} \left\{-H_{n}^{\prime}(\theta) \right\} = {\rm Var} \left\{
H_{n}(
\theta)\right\}$, holds under the given moment assumptions.
Therefore, the variance covariance matrix of $(\that,\hat{\sigma})$ reduces
to 
\[ \left( {\rm E }\left[ - H_{n}^{\prime}(\theta)\right]\right)^{-1} =
\left [ \begin{array}{lll}
       \left( 2 + \sigma^{-2} \right) \sm \bJ
\bJ^{T} &
\frac{2}{\sigma} \sm \bJ \\
\frac{2}{\sigma} \sm \bJ^{T}  & \frac{2 n}{\sigma^{2}}
\end{array}
\right]^{-1}. \]
The variance covariance matrix of $\that$, namely $V(\that)$, is given by the upper left corner of this matrix inverse.  Use standard
formulas for the inverse of a partitioned matrix to deduce~(\ref{var.that.dist}).

\end{document}